\documentclass[a4paper,11pt,leqno]{article}
\usepackage{amsmath}
\usepackage{amsfonts}
\usepackage{eucal}
\usepackage{amsthm}
\numberwithin{equation}{section}

\newcommand{\bigpare}[1]{\bigl(#1\bigr)}
\newcommand{\biggpare}[1]{\biggl(#1\biggr)}
\newcommand{\Bigpare}[1]{\Bigl(#1\Bigr)}

\newcommand{\bigbra}[1]{\bigl\{#1\bigr\}}

\newcommand{\biggbra}[1]{\biggl\{#1\biggr\}}

\newcommand{\Bigbrac}[1]{\Bigl[#1\Bigr]}
\newcommand{\biggbrac}[1]{\biggl[#1\biggr]}

\newcommand{\bigset}[2]{\bigl\{#1\bigm|#2\bigr\}}

\newcommand{\norm}[1]{\| #1 \|}
\newcommand{\bignorm}[1]{\bigl\| #1 \bigr\|}

\newcommand{\biggnorm}[1]{\biggl\| #1 \biggr\|}

\newcommand{\bigabs}[1]{\bigl| #1 \bigr|}
\newcommand{\Bigabs}[1]{\Bigl| #1 \Bigr|}
\newcommand{\biggabs}[1]{\biggl| #1 \biggr|}

\newcommand{\jap}[1]{\langle #1 \rangle}
\newcommand{\bigjap}[1]{\bigl\langle #1 \bigr\rangle}

\def\a{\alpha}
\def\b{\beta}
\def\c{\gamma}
\def\d{\delta}
\def\e{\varepsilon}

\def\g{\psi}

\def\l{\lambda}
\def\m{\mu}
\def\n{\nu}

\def\s{\sigma}

\def\x{\xi}
\def\y{\eta}
\def\z{\zeta}

\renewcommand{\L}{\Lambda}
\renewcommand{\O}{\Omega}

\newcommand{\OPS}{O\!P\!S}

\def\re{\mathbb{R}}

\def\ze{\mathbb{Z}}

\def\pa{\partial}

\newcommand{\supp}{\text{{\rm supp}\;}}

\newtheorem{thm}{Theorem}[section]
\newtheorem{lem}[thm]{Lemma}
\newtheorem{prop}[thm]{Proposition}
\newtheorem{cor}[thm]{Corollary}

\theoremstyle{definition}
\newtheorem{defn}{Definition}
\newtheorem{ass}{Assumption}

\theoremstyle{remark}
\newtheorem{rem}[thm]{Remark}

\newcommand{\WF}{W\!F}

\title{Semiclassical Singularity Propagation Property for 
Schr\"odinger Equations}
\author{Shu Nakamura%
\footnote{Graduate School of Mathematical Science, 
University of Tokyo, 3-8-1 Komaba, Meguro Tokyo, 
153-8914 Japan. 
E-mail:{\tt shu@ms.u-tokyo.ac.jp}.  Partially supported by JSPS Grant (B) 17340033.} }

\begin{document}
\maketitle

\begin{abstract}
We consider Schr\"odinger equations with variable coefficients, and it is 
supposed to be a long-range type perturbation of the flat Laplacian on $\re^n$. 
We characterize the wave front set of solutions to Schr\"odinger equations
in terms of the initial state. 
Then it is shown that the singularity propagates following the classical flow, 
and it is formulated in a semiclassical setting. 
Methods analogous to the 
long-range scattering theory, in particular a modified free propagator,  
are employed. 
\end{abstract}


\section{Introduction}
Let $H$ be a Schr\"odinger operator with variable coefficients:
\[
H= -\frac12 \sum_{j,k=1}^n \pa_{x_j} a_{jk}(x) \pa_{x_k} +V(x) 
\quad \text{on $L^2(\re^n)$}, 
\]
where $n\geq 1$ is the space dimension. Throughout this paper, we always assume
$a_{jk}(x)$ and $V(x)$ are real-valued $C^\infty$-class functions. 
Moreover, we assume:

\begin{ass}
For each $x\in\re^n$, $(a_{jk}(x))_{j,k}$ is a positive symmetric matrix. 
There is $\m>0$ such that for any multi-index $\a\in\ze_+^n$, there is 
$C_\a$ such that
\begin{align*}
\bigabs{\pa_x^\a\bigpare{a_{jk}(x)-\d_{jk}}} &\leq C_\a \jap{x}^{-\m-|\a|},  
\quad x\in\re^n, \\
\bigabs{\pa_x^\a V(x)} &\leq C_\a \jap{x}^{2-\m-|\a|},
\quad x\in\re^n. 
\end{align*}
\end{ass}

Then it is well-known that $H$ is essentially self-adjoint on $C_0^\infty(\re^n)$, 
and we denote the unique self-adjoint extension by the same symbol $H$. 
We let $u(t)=e^{-itH}u_0$ be the solution to the time-dependent 
Schr\"odinger equation:
\[
i\frac{\pa}{\pa t} u(t) = H u(t), \quad u(0)=u_0, \quad u_0\in L^2(\re^n).
\]
We study the microlocal singularity of $u(t)$. In particular, we characterize 
the wave front set of $u(t)$ in the nontrapping region, in terms of $u_0$. 
In order to describe our main result, we introduce several notations of the 
classical flow corresponding to $H$. Let $k(x,\x)$ be the classical kinetic 
energy, and let $p(x,\x)$ be the full Hamiltonian (modulo lower order terms):
\[
k(x,\x)=\frac12\sum_{j,k=1}^n a_{jk}(x)\x_j\x_k, \quad 
p(x,\x)=k(x,\x)+V(x), \quad 
x,\x\in\re^n.
\]
Let $\exp tH_p$ denote the Hamilton flow generated by a symbol $p$, i.e., 
if $(x(t),\x(t))=\exp tH_p (x_0,\x_0)$, then $(x(t),\x(t))$ is the solution
to the Hamilton equation:
\[
\frac{d}{dx} x(t) =\frac{\pa p}{\pa \x}(x(t),\x(t)), \quad
\frac{d}{dx} \x(t) =-\frac{\pa p}{\pa x}(x(t),\x(t)), \quad
t\in\re
\]
with $x(0)=x_0$, $\x(0)=\x_0$. 
\begin{defn}
For $(x_0,\x_0)\in\re^n\times\re^n$,we denote 
$(\tilde y(t),\tilde \y(t))=\exp tH_k(x_0,\x_0)$. 
$(x_0,\x_0)$ is called {\em backward nontrapping}\/ if 
$|\tilde y(t)|\to \infty$ as $t\to-\infty$. 
\end{defn}

For $a\in C^\infty(\re^{2n})$, we denote the Weyl quantization by
$a(x,D_x)$:
\[
a(x,D_x) u(x) =(2\pi)^{-n} \int e^{i(x-y)\cdot\x} a((x+y)/2,\x) u(y) \,dy\,d\x,
\]
where $u\in\mathcal{S}(\re^n)$ (see, e.g., H\"ormander \cite{Ho2}). 
We recall $(x_0,\x_0)\notin \WF(u)$, the wave front 
set of $u$, if and only if there exists $a\in C_0^\infty(\re^{2n})$ such that 
$a(x_0,\x_0)\neq 0$ and 
\[
\bignorm{a_h(x,D_x) u} = O(h^\infty)\quad \text{as $h\to 0$}, 
\]
where $a_h(x,\x)= a(x,h\x)$ (see, e.g., Martinez \cite{M}, Dimassi, Sj\"ostand \cite{DS}). 

\begin{thm}
Suppose $H$ satisfies Assumption~A, and let $u(t)=e^{-itH}u_0$, $u_0\in L^2(\re^n)$. 
Suppose, moreover, $(x_0,\x_0)$ is backward nontrapping, and let $t_0>0$.
Then $(x_0,\x_0)\notin \WF(u(t_0))$ if and only if there exists 
$a\in C_0^\infty(\re^{2n})$ such that $a(x_0,\x_0)\neq 0$ and 
\[
\bignorm{ \bigpare{ a_h\circ \exp t_0 H_p}(x,D_x) u_0} =O(h^\infty)
\quad \text{as $h\to 0$}. 
\]
\end{thm}

The main idea of the proof is very simple. Let $a\in C_0^\infty(\re^{2n})$ such that 
$a(x_0,\x_0)\neq 0$ and supported in a small neighborhood of $(x_0,\x_0)$. 
We note
\[
\bignorm{a_h(x,D_x)u(t_0)} = \bignorm{e^{it_0H} a_h(x,D_x) e^{-it_0 H} u_0}.
\]
If we formally apply the semiclassical Egorov theorem, we learn that the 
principal symbol of $e^{it_0H} a_h(x,D_x) e^{-it_0 H} u_0$  is given by 
$a_h\circ \exp t_0 H_p$, and we can obtain an asymptotic expansion of the 
symbol, where all the terms are supported in $\exp (-t_0 H_p)(\supp a)$. 
If this argument is justified, Theorem~1 follows immediately. 
However, in order to justify this argument in this framework, 
we need to find a suitable symbol class, 
which might be time-dependent. Instead of introducing time-dependent 
symbol class, we employ a scattering theoretical technique, which is an 
extension of the method used in \cite{N2}. 

Let $W(t,\x)$ be a solution to the momentum space Hamilton-Jacobi equation, 
which is constructed in Section~2. We study 
\[
\Omega(t) := e^{iW(t,D_x)} e^{-itH}
\]
instead of $e^{-itH}$ itself. Let 
\[
(y(t;x_0,\x_0),\y(t;x_0,\x_0)) = \exp tH_p (x_0,\x_0).
\]
If $(x_0,\x_0)$ is backward nontrapping, then it is shown in Section~2 that 
\begin{align*}
\x_-(-t_0,x_0,\x_0) 
&:= \lim_{\l\to+\infty} \l^{-1} \y(-t_0;x_0,\l \x_0),\\
z_-(-t_0;x_0,\x_0) 
&:= \lim_{\l\to+\infty} 
\biggpare{y(-t_0;x_0,\l \x_0)-\frac{\pa W}{\pa \x}(-t_0,\y(-t_0;x_0,\l\x_0))}
\end{align*}
exist. We will see that actually $\x_-$ and $z_-$ are  independent of $t_0$. 
We will show:

\begin{thm}
Suppose $H$ satisfies Assumption~A, and let $u(t)$, $(x_0,\x_0)$, $t_0>0$ as 
in Theorem~1. Then $(x_0,\x_0)\in\WF(u(t_0))$ if and only if 
\[
(z_-(-t_0;x_0,\x_0),\x_-(-t_0;x_0,\x_0))\in \WF\bigpare{e^{iW(-t_0,D_x)}u_0}.
\]
\end{thm}

Since the symbol of 
\[
e^{-iW(-t_0,D_x)} \bigpare{a_h\circ \exp t_0H_p}(x,D_x) e^{iW(-t_0,D_x)}
\]
is essentially supported in a small neighborhood of $(z_-,\x_-)$, 
Theorem~1.1 follows from Theorem~1.2 (see Subsection~3.4 for the detail). 
Theorem~1.2 is proved using an Egorov theorem for 
$\Omega(t) a_h(x,D_x) \Omega(t)^{-1}$. We note, at least formally, 
\begin{align*}
\frac{d}{dt} \Omega(t) &= i\frac{\pa W}{\pa t}(t,D_x) \Omega(t) 
- ie^{iW(t,D_x)} H e^{-itH} \\
&= -i \biggbra{ e^{iW(t,D_x)} H e^{-iW(t,D_x)} 
- \frac{\pa W}{\pa t}(t,D_x)} \Omega(t) \\
&=: -i L(t) \Omega(t). 
\end{align*}
Namely, $\Omega(t)$ is the evolution operator generated by a time-dependent 
self-adjoint operator $L(t)$. The principal symbol of $L(t)$ is given by 
\[
p\biggpare{x+\frac{\pa W}{\pa \x}(t,\x),\x} -\frac{\pa W}{\pa t}(t,\x) 
=p\biggpare{x+\frac{\pa W}{\pa \x}(t,\x),\x}- p\biggpare{\frac{\pa W}{\pa \x}(t,\x),\x}
\]
by virtue of the Hamilton-Jacobi equation.
This symbol is $O(\jap{\x}^{1-\m})$ if $t\neq 0$, and hence the speed of the 
propagation of singularity for $L(t)$ is 0 (away from $t=0$). 
However, at $t=0$, $L(0)$ has infinite propagation speed, 
and we observe a jump of the singularity. 
This propagation of singularity is described by the flow: 
$t \mapsto (z_-(t;x_0,\x_0), \x_-(t;x_0,\x_0))$, and we can conclude Theorem~1.2. 

Study of microlocal singularity of solutions to Schr\"odinger equation goes back 
at least to a work by Boutet de Monvel \cite{BM} (see also Lascar \cite{L}, 
Yamazaki \cite{Ym}, Zelditch \cite{Z}). Investigation to characterize the wave front 
set of $u(t)$ in terms of the initial state $u_0$ for variable coefficients 
Schr\"odinger equation was started by a work of Craig, Kappeler and Strauss \cite{CKS}. 
They showed that the microlocal regularity of the solution along a nontrapping 
geodesic follows from rapid decay of the initial state in a conic neighborhood of 
$-\x_-=-\lim_{t\to-\infty} \x(t)$. This property is called the microlocal smoothing 
property, and it was generalized and refined by Wunsch \cite{W}, Nakamura \cite{N1}
and Ito \cite{I}. Microlocal smoothing property in the analytic category was studied 
by Robbiano and Zuily \cite{RZ1, RZ2} and Martinez, Nakamura and Sordoni \cite{MNS}. 
Results in this paper may be considered as a refinement of these works, and the 
microlocal smoothing property (in the $C^\infty$-category) follows immediately 
from Theorem~1.1. Similar characterization of wave front set for solutions to Schr\"odinger 
equation is recently obtained by Hassel and Wunsch \cite{HW}. They considered the problem 
in the framework of {\it scattering metric}\/, and the assumptions and the 
proof are quite different. In a previous paper, 
Nakamura \cite{N2} considered the case of short-range 
perturbations, i.e., $\m>1$, and the results in this paper are its generalizations.

On the other hand, the singularity of solutions to perturbed 
harmonic oscillator Schr\"odinger 
equation was studied by Zelditch \cite{Z}, Yajima \cite{Yj},  Kapitanski, Rodnianski
and Yajima \cite{KRY} and Doi \cite{D1,D2}. 
The idea of these papers, especially those by Doi, is closely related to our proof. 

Recently, Strichartz estimates for variable coefficient Schr\"odinger operator 
was studied by several authors, e.g., Staffilani and Tataru \cite{ST}, 
Robbiano and Zuily \cite{RZ3}, Burq, G\'erard and Tzvetkov \cite{BGT}, 
Bouclet and Tzvetkov \cite{BT}. Strichartz estimate is another expression of the 
smoothing property of Schr\"odinger equations, and there should be implicit 
relationship with our results. In particular, Bouclet and Tzvetkov used 
the Isozaki-Kitada modifier to treat long-range perturbations, and it is 
analogous to our modified free propagator, though the formulation and the 
construction are completely different. 

The paper is organized as follows: In Section~2, we consider the classical 
motions generated by the kinetic energy and the total Hamiltonian. 
In particular, we construct a solution to the momentum space Hamilton-Jacobi 
equation and show the existence of the modified classical wave operator. 
We prove Theorems~1.2 and then Theorem~1.1 in Section~3. 

Throughout this paper, we use the following notation: $S(m,g)$ denotes the 
H\"ormander symbol class (cf. H\"ormander~\cite{Ho2}, Chapter~18). 
For a compact set $K\subset \re^n$, $S_K(m,g)$ denotes the same symbol class 
restricted to functions on $K\times\re^n$. For a symbol $a(x,\x)$, $a(x,D_x)$ 
denotes the Weyl quantization of $a$. 

\noindent
\textbf{Acknowledgment.} The author would like to thank the referee for pointing out 
numerous errors in the first version (his apology), and also for providing valuable suggestions. 
He also thanks Kenji Yajima, Andr\'e Martinez, Shin-ichi Doi and Ken-ichi Ito for valuable discussions. 


\section{Hamilton flows and solution to the Hamilton-Jacobi equation}

\subsection{Properties of nontrapping geodesic flow}

Here we consider the Hamilton flow for the kinetic energy: 
$k(x,\x) =\frac12\sum a_{jk}(x)\x_j \x_k$. 
We always suppose Assumption~A is satisfied. 

\begin{prop}
Let $(x_0,\x_0)\in\re^{2n}$ and suppose $(x_0,\x_0)$ is backward nontrapping. 
Then there exists $C>0$ such that 
\[
|\tilde y(t)| \geq C^{-1}|t|-C, \qquad t\leq 0,
\]
where 
\[
(\tilde y(t), \tilde\y(t)) =(\tilde y(t;x_0,\x_0),\tilde \y(t;x_0,\x_0)) = 
\exp tH_k (x_0,\x_0).
\]
Moreover, $C$ is taken locally uniformly with respect to $(x_0,\x_0)$
\end{prop}

\begin{proof}
At first we recall the conservation of the energy:
\[
k(\tilde y(t),\tilde \y(t)) 
=\frac12 \sum_{j,k} a_{jk}(\tilde y(t))\tilde \y_j(t) \tilde \y_k(t)
=k(x_0,\x_0).
\]
By the uniform ellipticity of $k(x,\x)$, we learn that there exists $C_1>0$ 
such that 
\[
C_1^{-1} \leq |\tilde\y(t)|\leq C_1, \qquad t\in\re.
\]
We compute 
\begin{align*}
\frac{d^2}{dt^2} |\tilde y(t)|^2 &= 2\frac{d}{dt}\biggpare{\tilde y(t)\cdot 
\frac{d \tilde y}{dt}(t)} 
=2\frac{d}{dt}\biggpare{\sum_{j,k} a_{jk}(\tilde y(t))\tilde y_j(t) \tilde \y_k(t)}\\
&= 4k(\tilde y(t),\tilde\y(t))+ \tilde U(\tilde y(t),\tilde\y(t)), 
\end{align*}
where
\begin{align*}
\tilde U(x,\x) &= 2\sum_{j,k,\ell} a_{jk}(x)\bigpare{a_{j\ell}(x)-\d_{j\ell}}\x_\ell\x_k \\
&\quad - \sum_{j,k,\ell,m}a_{jk}(x) \frac{\pa a_{\ell m}}{\pa x_k}(x)x_j\x_\ell\x_m 
+2\sum_{j,k,\ell,m} \frac{\pa a_{jk}}{\pa x_\ell}(x) a_{\ell m}(x)x_j \x_k \x_m.
\end{align*}
By Assumption~A, it is easy to see 
\[
|\tilde U(x,\x)|\leq C\jap{x}^{-\m}|\x|^2, 
\]
and this implies
\[
\frac{d^2}{dt^2}|\tilde y(t)|^2 \geq 4k(x_0,\x_0)
-C\jap{\tilde y(t)}^{-\m}|\tilde\y(t)|^2.
\]
We can choose $R>0$ so large that 
\[
4k(x_0,\x_0)-C R^{-\m} C_1^2 \geq \e>0.
\]
Since $(x_0,\x_0)$ is backward nontrapping, there exists $t_0<0$ such that 
\[
|\tilde y(t_0)|=R, \quad \frac{d}{dt}|\tilde y(t_0)|\leq 0,
\]
and hence 
\[
\frac{d^2}{dt^2}|\tilde y(t_0)|^2 \geq \e.
\]
Then by the convexity of $|\tilde y(t)|^2$, we conclude 
\[
|\tilde y(t)|^2 \geq R+ \e (t_0-t)^2/2, \qquad t\leq t_0, 
\]
and the assertion follows immediately. 
\end{proof}

\begin{prop}
Suppose $(x_0,\x_0)$ is backward nontrapping. Then 
\[
\x_- := \lim_{t\to-\infty} \tilde \y(t;x_0,\x_0)
\]
exists. 
\end{prop}

\begin{proof}
By Proposition~2.1 and Assumption~A, we learn 
\begin{align*}
\frac{d}{dt} \tilde \y_j(t) &=  -\frac12 \sum_{k,\ell} (\pa_{x_j}a_{k\ell})
(\tilde y(t)) \tilde \y_k(t)\tilde \y_\ell(t) \\
&= O(|\tilde y(t)|^{-1-\m}) = O(\jap{t}^{-1-\m})
\end{align*}
as $t\to-\infty$. Hence 
\[
\x_-= \lim_{t\to-\infty} \tilde \y(t) = \x_0-\int_{-\infty}^0 \frac{d}{dt}\tilde \y(t) dt
\]
exists. 
\end{proof}

By the above proof, we also observe 
\[
|\x_--\tilde\y(t)| \leq C\jap{t}^{-\m}, \quad t\to-\infty,
\]
and $C$ can be taken locally uniformly in $(x_0,\x_0)$. 
Let $0<\d_1<1$, $R>0$, and we set
\[
\O_{R,\d_1}=\bigset{(x,\x)\in\re^{2n}}{R-1<|x|<R+1,\tfrac12<|\x|<2, 
x\cdot\x\leq -\d_1|x|\;|\x|}
\]
be a neighborhood of $\bigset{(x,-x/|x|)\in\re^{2n}}{|x|=R}$. 
We fix $\d_1>0$. If $R$ is sufficiently large, we have 
\begin{align*}
\frac{d}{dt}|\tilde y(t)|^2 \;\Big|_{t=0} 
&= \sum_{j,k} a_{jk}(x)x_j\x_k \\
&=x\cdot\x +\sum_{j,k} (a_{jk}(x)-\d_{jk})x_j\x_k \\
&\leq -\d_1|x|\;|\x|+ \frac{\d_1}{2}|x|\;|\x| 
=-\frac{\d_1}{2}|x|\;|\x|
\end{align*}
for $(x,\x)\in\Omega_{R,\d_1}$. Hence, in particular, $(x,\x)$ is 
backward nontrapping and 
\[
|\tilde y(t)|^2 \geq |x|^2 +\frac{\d_1}{8}|x|\;|t| +\e|t|^2 
\qquad \text{for $t\leq 0$}.
\]
Thus we have proved the following assertion:

\begin{prop}
Let $0<\d_1<1$. There exist $R_0>0$ and $\d_2>0$ such that 
if $R\geq R_0$ then 
\[
|\tilde y(t;x,\x)| \geq |x|+\d_2 |t|, \qquad t\leq 0,\;\;(x,\x)\in \Omega_{R,\d_1}. 
\]
\end{prop}

We note that since $k(x,\x)$ is homogeneous in $\x$, the flow also has the following 
homogeneity: for $\l>0$, 
\begin{align*}
\tilde y(t;x,\l\x) &= \tilde y(\l t;x,\x), \\
\tilde\y(t;x,\l \x) &= \l \tilde\y(\l t;x,\x).
\end{align*}
Thus we learn the following property concerning the high energy 
asymptotics of the geodesic flow:

\begin{prop}
(i) Suppose $(x_0,\x_0)$ is backward nontrapping. Then for any $t<0$, 
$\l>0$, 
\[
|\tilde  y(t;x_0,\l\x_0)|\geq C^{-1}\l |t|-C, 
\]
and 
\[
\x_-(x_0,\x_0) =\lim_{\l\to+\infty} \l^{-1} \tilde \y(t;x_0,\l \x_0)
\]
exists. $\x_-$ is independent of $t<0$.
\newline
(ii) Let $0<\d_1<1$. Then there exist $R_0>0$ and $\d_2>0$ such that if 
$R\geq R_0$ then
\[
|\tilde y(t;x,\x)|\geq |x|+\d_2|t|\;|\x|
\]
for $t\leq 0$ and  
\[
(x,\x) \in\bigset{(x,\x)\in\re^{2n}}{R-1<|x|<R+1, x\cdot\x \leq -\d_1|x|\;|\x|}.
\]
In particular, 
\[
\x_-(x,\x) =\lim_{\l\to+\infty} \l^{-1}\tilde \y(t;x,\l\x), 
\qquad (x,\x)\in \Omega_{R,\d_1},
\]
converges uniformly in $\Omega_{R,\d_1}$. 
\end{prop}


\subsection{High energy asymptotics of the Hamilton flow}

Now we consider the Hamilton flow: 
\[
(y(t;x,\x),\y(t;x,\x)) =\exp tH_p (x,\x).
\]
We recall $(y(t),\y(t))$ satisfies the Hamilton equation:
\begin{align*}
\frac{d}{dt} y_j(t)&= \sum_{k=1}^n a_{jk}(y(t)) \,\y_k(t), \\
\frac{d}{dt} \y_j(t) 
&= -\frac12 \sum_{j,k=1}^n \frac{\pa a_{k\ell}}{\pa x_j}(y(t))\, \y_k(t)\,\y_\ell(t) 
-\frac{\pa V}{\pa x_j}(y(t)).
\end{align*}
At first we prepare an a priori estimate:

\begin{prop}
Let $T>0$. Then there exist $\a, \b,\c>0$ such that 
\[
|y(t;x,\x)| \leq \a |\x|, \quad |\y(t;x,\x)|\leq \b|\x|
\]
if $|\x|>1$, $t\in [-T,T]$ and $|x|\leq \c|\x|$.
\end{prop}

\begin{proof}
We note
\[
p(x,\x)= k(x,\x)+V(x) \leq c_1\jap{\x}^2
\]
with some $c_1>0$ if $|x|\leq \c|\x|$. The by the conservation of energy, we learn 
\begin{align*}
|\y(t;x,\x)|&\leq c_2 \sqrt{k(y,\y)} =c_2 \sqrt{p(y,\y)-V(y)} \\
&\leq c_3 (\jap{\x}+\jap{y}^{((2-\m)/2}) \leq c_3 (\jap{\x}+\jap{y}).
\end{align*}
Hence we have 
\[
\Bigl| \dfrac{d}{dt} y(t;x,\x)\Bigr| \leq c_3 (\jap{\x}+\jap{y}).
\]
By using the Duhamel formula, we obtain
\[
|y(t)|\leq e^{c_3t}|x|+ \int_0^t e^{c_3(t-s)}c_3 \jap{\x} ds 
\leq c_4 \jap{\x}
\]
if $0\leq t\leq T$ and $|x|\leq \c |\x|$. Then under the same assumption we also have 
\[
|\y(t)|\leq c_3(\jap{\x}+\jap{y})\leq c_5 \jap{\x}.
\]
The case $-T\leq t\leq 0$ is similar, and we omit the detail. 
\end{proof}

If we denote 
\begin{align*}
y^\l(t;x,\x) &= y(t/\l;x,\l \x), \\
\y^\l(t;x,\x) &= \frac1\l \,\y(t/\l;x,\l\x),
\end{align*}
for $\l>0$, then $(y^\l(t),\y^\l(t))$ satisfies
\begin{align*}
\frac{d}{dt} y^\l_j(t) &= \sum_k a_{jk}(y^\l)\,\y^\l_k, \\
\frac{d}{dt}\y^\l_j(t) &= -\frac12\sum_{k,\ell} \frac{\pa a_{k\ell}}{\pa x_j}(y^\l)
\, \y^\l_k\,\y^\l_\ell -\frac{1}{\l^2}\frac{\pa V}{\pa x_j}(y^\l), 
\end{align*}
with the initial condition: $y^\l(0)=x$, $\y^\l(0)=\x$. 
By the continuity of the solutions to ODE's in the coefficients, 
we learn 
\[
y^\l(t)\to \tilde y(t), \quad \y^\l(t) \to \tilde \y(t) \quad \text{as }\l\to+\infty,
\]
locally uniformly in $t\in\re$. In particular, if $(x,\x)$ is nontrapping, 
then for any $R>0$, $|y^\l(t)|>R$ for $t\ll 0$ and $\l\gg 0$. In fact, we have the 
following stronger assertion:

\begin{prop}
Suppose $(x,\x)$ is backward nontrapping, and let $t_0<0$. 
Then there exist $C>0$ and $\l_0>0$ such that 
\[
|y^\l(t)|\geq C^{-1} |t| -C, \qquad \text{for}\quad
\l t_0\leq t\leq 0,\;  \l\geq \l_0, 
\]
where $y^\l(t)=y^\l(t;x,\x)$. Moreover, $C$ can be taken locally uniformly 
with respect to $(x,\x)$. 
\end{prop}

\begin{proof}
The proof is analogous to Proposition~2.1. 
By Proposition~2.5, we have 
\[
|y^\l(t)|\leq \a\, \l|\x|
\quad \text{for } \l t_0\leq t\leq 0, 
\]
if $\l$ is sufficiently large (so that $|x|\leq \b\,\l|\x|$). 
As in the proof of Proposition~2.1, we have 
\[
\frac{d^2}{dt^2} |y^\l(t)|^2 =4p^\l(y^\l(t),\y^\l(t)) +U(y^\l(t),\y^\l(t)), 
\]
where 
\begin{align*}
p^\l(x,\x)&= \frac12\sum_{j,k=1}^n a_{jk}(x)\,\x_j\,\x_k +\frac{1}{\l^2} V(x), \\
U(x,\x) &=\tilde U(x,\x) -\frac{4}{\l^2}V(x) 
-\frac{2}{\l^2}\sum_{j,k} a_{jk}(x)\,x_j\,\frac{\pa V}{\pa x_k}(x). 
\end{align*}
These imply 
\[
\frac{d^2}{dt^2} |y^\l(t)|^2 \geq 4k(x,\x) -C\l^{-\m} -C\jap{y^\l(t)}^{-\m}. 
\]
Then, by noting the above remark that $y^\l(t)\to \tilde y(t)$ as $\l\to+\infty$, 
the same argument as in the proof of Proposition~2.1 applies, 
and we conclude the assertion.
\end{proof}

\begin{cor}
Let $(x,\x)$, $t_0$, $C$ and $\l_0$ as in Proposition~2.6. Then 
\[
|y(t;x,\l\x)|\geq C^{-1}\,\l|t|-C 
\quad \text{for } t_0\leq t\leq 0, \;\l\geq \l_0.
\]
\end{cor}

As well as Proposition~2.4, we also have the following proposition:

\begin{prop}
Let $0<\d_1<1$ and $t_0<0$. Then there exist $R_0>0$, $\d_2>0$ and $\l_0>0$ 
such that if $R\geq R_0$ then 
\[
|y(t;x,\x)|\geq |x|+\d_2\,|t|\,|\x|, \qquad t_0\leq t\leq 0, 
\]
for $(x,\x)\in\bigset{(x,\x)}{R-1<|x|<R+1,|\x|\geq \l_0,x\cdot\x\leq -\d_1|x|\,|\x|}$. 
\end{prop}

\begin{prop}
Suppose $(x,\x)$ is backward nontrapping. Then for any $t_0<0$, there exists $C>0$ such that 
\begin{align*}
|\y(t;x,\l\x)-\tilde\y(t;x,\l\x)|&\leq C\l^{1-\m}|t|^{2-\m}, \\
|y(t;x,\l\x)-\tilde y(t;x,\l\x)| &\leq C \l^{1-\m} |t|^{3-\m}
\end{align*}
for $t\in [t_0,-1/\l]$ and $\l>1$. 
\end{prop}

\begin{proof}
It suffices to show the equivalent assertion:
\begin{align*}
|\y^\l(t;x,\x)-\tilde\y(t;x,\x)| &\leq C\l^{-2}|t|^{2-\m}, \\
|y^\l(t;x,\x)-\tilde y(t;x,\x)| &\leq C\l^{-2}|t|^{3-\m}
\end{align*}
for $t\in [\l t_0,-1]$. 
By the Hamilton equation, we have 
\begin{align*}
&\frac{d}{dt}\Bigpare{\y^\l_j(t)-\tilde\y_j} = -\frac12 \sum_{k,\ell}
\biggpare{\frac{\pa a_{k\ell}}{\pa x_j}(y^\l)\,\y_k^\l\,\y^\l_\ell 
-\frac{\pa a_{k\ell}}{\pa x_j}(\tilde y)\,\tilde \y_k\,\tilde \y_\ell}
-\frac{1}{\l^2}\frac{\pa V}{\pa x_j}(y^\l), \\
&\frac{d}{dt}\Bigpare{y^\l_j(t)-\tilde y_j(t)} = \sum_k  \Bigpare{a_{jk}(y^\l)\,\y^\l_k 
-a_{jk}(\tilde y)\,\tilde \y_k}. 
\end{align*}
These imply
\begin{align*}
\biggabs{\frac{d}{dt}\Bigpare{\y^\l-\tilde\y}} &\leq c_1 \Bigpare{
|t|^{-1-\m}|\y^\l-\tilde\y|+|t|^{-2-\m}|y^\l-\tilde y|+\l^{-2}|t|^{1-\m}}, \\
\biggabs{\frac{d}{dt}\Bigpare{y^\l-\tilde y}} 
&\leq c_1 |\y^\l-\tilde \y|+c_1|t|^{-1-\m}|y^\l-\tilde y|, 
\end{align*}
for $t\leq -1$ with some $c_1>0$ (cf.\ Lemma~A.1 in Appendix). If $t\leq -T<0$, we have 
\begin{align*}
\biggabs{\frac{d}{dt}\Bigpare{\y^\l-\tilde\y}} &\leq c_1 \Bigpare{
T^{-\m} |t|^{-1}|\y^\l-\tilde\y|+T^{-\m} |t|^{-2}|y^\l-\tilde y|+\l^{-2}|t|^{1-\m}}, \\
\biggabs{\frac{d}{dt}\Bigpare{y^\l-\tilde y}} 
&\leq c_1 |\y^\l-\tilde \y|+c_1T^{-\m} |t|^{-1}|y^\l-\tilde y|. 
\end{align*}
Thus, for $t<-T$, $|\y^\l-\tilde\y|$ and $|y^\l-\tilde y|$ are majorized by a solution to 
\begin{align*}
-Z' &\geq c_1\Bigpare{T^{-\m} |t|^{-1}Z + T^{-\m} |t|^{-2} Y + \l^{-2}|t|^{1-\m}}, \\
-Y' &\geq Z+ c_1 T^{-\m} |t|^{-1} Y
\end{align*}
with 
\[
Z(-T) \geq |\y^\l(-T) -\tilde\y(-T)|, \quad 
Y(-T)\geq |y^\l(-T)-\tilde y(-T)|. 
\]
If we set
\[
Y(t) =c_2\l^{-2} |t|^{3-\m}, \quad Z(t)=c_3\l^{-2}|t|^{2-\m}, 
\]
then the differential inequalities are satisfied if 
\begin{align*}
c_3(2-\m) &\geq c_1\bigpare{c_3T^{-\m}+c_2T^{-\m} +1}; \\
c_2(3-\m) &\geq c_1 c_3 +c_1c_2 T^{-\m}.
\end{align*}
In other words, if 
\begin{align*}
c_3((2-\m)-c_1T^{-\m}) &\geq c_1c_2T^{-\m}+c_1; \\
c_2((3-\m)-c_1T^{-\m}) &\geq c_1 c_3. 
\end{align*}
We choose $T$ so large that 
\[
((2-\m)-c_1T^{-\m})^{-1}\times c_1^2 T^{-\m}((3-\m)-c_1T^{-\m})^{-1} <1, 
\]
and set $c_3=((3-\m)-c_1T^{-\m})c_1^{-1} c_2$. If $c_2$ is sufficiently large, 
the above inequalities are satisfied. 

Since $|y^\l(-T)-\tilde y(-T)|$, $|\y^\l(-T)-\tilde\y(-T)|=O(\l^{-2})$ as 
$\l\to+\infty$, the initial condition is also satisfied if $c_2$ is taken 
sufficiently large. Thus we conclude the assertion for $t\in [\l t_0,-T]$. 
The estimate for $t\in [-T,-1]$ is obvious. 
\end{proof}

Proposition~2.9 implies, in particular, 
\[
\lim_{\l\to+\infty}\l^{-1}\y(t;x,\l\x) 
=\lim_{\l\to+\infty} \l^{-1}\tilde\y(t;x,\l\x) =\x_-(x,\x). 
\]


\subsection{Construction of a solution to the Hamilton-Jacobi equation}

In order to construct a solution to the momentum space Hamilton-Jacobi equation, 
we prepare one more lemma about the classical flow:

\begin{prop}
Let $\d_1>0$ and $t_0<0$. There exist $R_0>0$, $c_0>0$  and $C>0$ such that 
\begin{align*}
\biggabs{\frac{\pa}{\pa x}\y(t;x,\x)}&\leq C R^{-1-\m}|\x|, \\
\biggabs{\frac{\pa}{\pa\x}(\y(t;x,\x)-\x)} &\leq C R^{-\m}
\end{align*}
for $t_0\leq t\leq 0$, 
\[
(x,\x)\in\Omega :=\bigset{(x,\x)\in\re^{2n}}{\bigabs{|x|-R}\leq 1, |\x|\geq \l, 
x\cdot\x\leq -\d_1|x|\cdot|\x| }
\]
with $R\geq R_0$ and $\l\geq c_0 R$. Moreover, for any $\a,\b\in\ze_+^n$, there is 
$C_{\a\b}>0$ such that 
\begin{align*}
\biggabs{\biggpare{\frac{\pa}{\pa x}}^\a\biggpare{\frac{\pa}{\pa\x}}^\b \bigpare{y(t;x,\x)-x}} 
&\leq C_{\a\b}\,|t|\,\jap{\x}^{1-|\b|}, \\
\biggabs{\biggpare{\frac{\pa}{\pa x}}^\a\biggpare{\frac{\pa}{\pa\x}}^\b \bigpare{\y(t;x,\x)-\x}} 
&\leq C_{\a\b}\,\jap{\x}^{1-|\b|},
\end{align*}
for $(x,\x)\in\Omega$ and $t\in [t_0,0]$. 
\end{prop}

\begin{proof}
We set $\l=|\x|$ and consider
\begin{align*}
y^\l(t;x,\x) &=y(t/\l;x,\l \x), \\
\y^\l(t;x,\x) &= \l^{-1} \y(t/\l;x,\l \x).
\end{align*}

Then it suffices to show the above estimates for $\y^\l$ and $y^\l$ with 
$|\x|=1$, $\l\geq \l_0$ and $t\in [\l t_0,0]$.

We mimic the argument of H\"ormander \cite{Ho} Lemma~3.7. Let $s$ be a variable 
$x_j$ or $\x_j$, $j=1,\dots,n$. By the Hamilton equation, we have 
\begin{align}
\frac{d}{dt}\biggpare{\frac{\pa y^\l_j}{\pa s}} &= \sum_{k,\ell} 
\frac{\pa a_{jk}}{\pa x_\ell}(y^\l)\frac{\pa y^\l_\ell}{\pa s} \,\y_k^\l
+\sum_k a_{jk}(y^\l) \frac{\pa \y^\l_k}{\pa s}, \\
\frac{d}{dt}\biggpare{\frac{\pa \y^\l_j}{\pa s}} &= -\frac12\sum_{k,\ell,m}
\frac{\pa^2 a_{k\ell}}{\pa x_j \pa x_m}(y^\l)\, \y^\l_k \, \y^\l_\ell \,
\frac{\pa y_m^\l}{\pa s} \\
& \quad -\sum_{k,\ell} \frac{\pa a_{k\ell}}{\pa x_j}(y^\l)\, \y^\l_k\,
\frac{\pa \y^\l_\ell}{\pa s} -\frac{1}{\l^2} 
\sum_k \frac{\pa^2 V}{\pa x_k \pa x_j}(y^\l) \frac{\pa y^\l_k}{\pa s}. \nonumber
\end{align} 
Then $|{\pa y^\l}/{\pa s}|$ and $|{\pa \y^\l}/{\pa s}|$ are 
majorized by a solution to 
\begin{align*}
-\frac{d}{dt} Y &\geq c_1 (R+\d|t|)^{-1-\m} Y +c_1 Z, \\
-\frac{d}{dt} Z &\geq c_1 (R+\d|t|)^{-2-\m} Y + c_1 (R+\d|t|)^{-1-\m} Z 
+ \frac{c_1}{\l^2} (R+\d|t|)^{-\m} Y,
\end{align*}
with $Y(0)\geq 0$ and $Z(0)\geq 1$ if $s=\x_j$, $Y(0)\geq 1$ and 
$Z(0)\geq 0$ if $s=x_j$. Note we consider the inequality in $t<0$. 

We set 
\[
Y=c_2(R-\d t), \quad Z=c_3(1-(R-\d t)^{-\m'}), \quad \l t_0\leq t\leq 0
\]
with $0<\m'<\m$. Then the differential inequalities for the majorants are satisfied if
\begin{align}
c_2\d &\geq c_1 c_2 R^{-\m} +c_1 c_3, \\
c_3 \d \m' &\geq R^{-(\m-\m')} \biggpare{c_1 c_2  + c_1 c_3 
+ c_1 c_2 \biggpare{\frac{ R-\d\l t_0}{\l}}^2},
\end{align}
and $R^{-\m'}\leq 1/2$ so that $Z>0$. 
We note 
\[
\biggabs{\frac{R-\d\l t_0}{\l}} = \biggabs{\frac{R}{\l} -\d t_0} \leq c_0^{-1} +\d |t_0|
\]
since $\l>c_0 R$. Now we choose $c_2/c_3=\c>2c_1/\d$, and choose $R_0$ so that 
\[
R_0 > \max \biggbra{ 2^{1/\m'}, \biggpare{\frac{2c_1}{\d}}^{1/\m}, 
\biggpare{\frac{\c c_1 c_2}{\d\m'}\Bigbrac{1+\c^{-1}+(c_0^{-1}+\d |t_0|)^2}}^{1/(\m-\m')}}, 
\]
then the above conditions are satisfied. Thus we learn 
\[
\biggabs{\frac{\pa y_j^\l}{\pa s}(t)} \leq c_2 (R-\d t), \qquad
\biggabs{\frac{\pa \y_j^\l}{\pa s}(t)} \leq c_2/\c, 
\]
for $R\geq R_0$, $\l\geq c_0 R$ and $t\in [\l t_0,0]$, provided 
\[
\biggabs{\frac{\pa y_j^\l}{\pa s}(0)} \leq c_2 R, 
\qquad \biggabs{\frac{\pa \y_j^\l}{\pa s}(0)} \leq c_2/(2\c).
\]
We now consider the case $s=x_k$. Then we may set $c_2=R^{-1}$ and we have 
\[
\biggabs{\frac{\pa y_j^\l}{\pa x_k}(t)} \leq 1-\frac{\d t}{R}, \qquad 
\biggabs{\frac{\pa \y_j^\l}{\pa x_k}(t)} \leq \frac{1}{\c R}.
\]
We integrate the equation (2.2) again to obtain 
\begin{align*}
\biggabs{\frac{\pa \y_j^\l}{\pa x_k}(t)} &\leq \frac{c_1}{R} \int_t^0 (R-\d r)^{-1-\m}dr
+ \frac{c_1}{\c R}\int_t^0 (R-\d r)^{-1-\m} dr \\
& \qquad +\frac{c_1}{R\l^2} \int_t^0 (R-\d r)^{1-\m} dr  \\
&\leq (\biggpare{\frac{c_1}{\m}+\frac{c_1}{\c\m}} R^{-1-\m} 
+ c_1 \biggpare{\frac{R-\d t}{\l}}^2 R^{-1-\m} \\
&\leq C R^{-1-\m}
\end{align*}
if $t\in [\l t_0,0]$ and $\l\geq c_0R$. Similarly, if $s=\x_k$, we may set 
$c_2=2\c$ and we have 
\[
\biggabs{\frac{\pa y_j^\l}{\pa \x_k}(t)} \leq 2\c (R-\d t), \qquad 
\biggabs{\frac{\pa \y_j^\l}{\pa \x_k}(t)} \leq 2.
\]
By integrating the equation (2.2), we conclude 
\[
\biggabs{\frac{\pa \y_j^\l}{\pa \x_k}(t) -\d_{jk}} \leq C' R^{-\m}.
\]

For higher derivatives, we prove the estimates by induction. It suffices to show
\begin{align*}
\bigabs{ \pa_x^\a \pa_\x^\b (y^\l(t;x,\x)-x)} &\leq C_{\a\b}|t|, \\
\bigabs{ \pa_x^\a \pa_\x^\b (\y^\l(t;x,\x)-\x)} &\leq C_{\a\b}
\end{align*}
for $t\in [\l t_0,0]$. We suppose these hold for $|\a+\b|<k$, and 
let 
\[
Y(t)={ \pa_x^\a \pa_\x^\b (y^\l(t;x,\x)-x)}, \quad 
Z(t) ={ \pa_x^\a \pa_\x^\b (\y^\l(t;x,\x)-\x)}
\]
with $|\a+\b|=k$. Then by the induction hypothesis, we can show $Y$ and $Z$ satisfy
\begin{align*}
& Y' = A_{11} Y+ A_{12} Z + A_{13}, \\
&Z' = A_{21} Y + A_{22} Z +A_{23} + \l^{-2} (A_{31}Y +A_{33}), \\
&Y(0)=Z(0)=0, 
\end{align*}
where 
\begin{align*}
&A_{11}=O(\jap{t}^{-1-\m}), \quad A_{12}=O(1), \quad A_{13}=O(\jap{t}^{-\m}), \\
&A_{21}=O(\jap{t}^{-2-\m}), \quad A_{22} = O(\jap{t}^{-1-\m}), \quad A_{23} = O(\jap{t}^{-1-\m}), \\
&A_{31}= O(\jap{t}^{-\m}), \quad A_{32} =O(\jap{t}^{1-\m}), 
\end{align*}
which itself is proved by induction. Then for $t\in [\l t_0,-1]$, we have 
\begin{align*}
\bigabs{Y'} &\leq c_1 (\jap{t}^{-1-\m}|Y| +|Z| + \jap{t}^{-\m}), \\
\bigabs{Z'} &\leq c_1 (\jap{t}^{-2-\m}|Y| +\jap{t}^{-1-\m} |Z| + \jap{t}^{-1-\m}). 
\end{align*}
These imply $Y$ and $Z$ are majorized by $M\jap{t}$ and $M$, respectively, with sufficiently large $M$
(the proof is analogous to the above argument). 
By integrating the differential equation again, we conclude the assertion for $|\a+\b|=k$. 
\end{proof}

We note the above proof for the derivatives works for $(\tilde y(t;x,\x), \tilde \y(t,x,\x))$  ($t<0$) 
if $(x,\x)$ is backward nontrapping. In particular, we learn that $\pa_t \pa_x^\a\pa_\x^\b \tilde\y(t;x,\x)$ 
is integrable with respect to $t$ in $(-\infty,0]$, and hence we conclude $\pa_x^\a\pa_\x^\b \tilde\y(t;x,\x)$ 
converges as $t\to-\infty$, and the estimate is locally uniform.  Thus we have 

\begin{cor} Suppose $(x,\x)$ is backward nontrapping. Then 
\[
(x,\x) \mapsto \x_-(x,\x)
\]
is a $C^\infty$ map, and $\tilde\y(t;x,\x)$ converges to $\x_-(x,\x)$ 
locally uniformly with all the derivatives as $t\to-\infty$. 
\end{cor}

Now we consider the map:
\[
\Lambda: \x\;\longmapsto\; \y(t;-R\x/|\x|,\x).
\]
Proposition~2.10 implies $\bignorm{\frac{\pa \Lambda}{\pa\x}-I}=O(R^{-\m})$ 
uniformly for $|\x|\geq c_0 R$. We choose $R$ so large that 
$\frac{\pa \Lambda}{\pa\x}$ is invertible for $|\x|\geq c_0R$. It is also easy to 
see that $|\Lambda -\x|=O(R^{-\m}|\x|)$ for $|\x|\geq c_0 R$, and hence 
$\mbox{Ran}\Lambda \supset \bigset{\x\in\re^n}{|\x|\geq c_4 R}$ with some $c_4>0$. 
Then we set 
\[
\z(t,\cdot) =\Lambda(t,\cdot)^{-1}\;:\; 
\bigset{\x}{|\x|\geq c_4 R} \longrightarrow \re^n, 
\]
i.e., 
\[
\y(t;-R\z(t,\x)/|\z(t,\x)|,\z(t,\x)) =\x \qquad \text{for } |\x|\geq c_4 R.
\]
By Proposition~2.10, we learn
\begin{equation}
\biggabs{\biggpare{\frac{\pa }{\pa \x}}^\a \z(t,\x)}\leq C_\a \jap{\x}^{1-|\a|}, 
\qquad t\in[t_0,0], \; |\x|\geq c_4 R. 
\end{equation}
Then we set 
\[
W_1(t,\x) =\int_0^t \bigpare{p(y(s),\y(s))+y(s)\cdot\pa_t \y(s)}ds -R|\x|, \qquad |\x|\geq c_4 R, 
\]
where 
\[
y(s) = y(s;-R\z(t,\x)/|\z(t,\x)|,\z(t;\x)), \quad  
\y(s) = \y(s;-R\z(t,\x)/|\z(t,\x)|,\z(t;\x)).
\]
It is well-known that $W_1(t,\x)$ satisfies the Hamilton-Jacobi 
equation (cf. Reed-Simon~\cite{RS} Section~XI.9):
\[
\frac{\pa}{\pa t} W_1(t,\xi) =p\biggpare{\frac{\pa W_1}{\pa \x}(t,\x),\x},
\qquad |\x|\geq c_4 R.
\]
By the construction we have 
\[
\pa_\x W_1(t,\x) =y(t;-R\z(t,\x)/|\z(t,\x)|,\z(t;\x)), 
\]
and 
\begin{equation}
\bigabs{\pa_\x^\a W_1(t,\x)} \leq C_\a \jap{\x}^{2-|\a|}, 
\qquad t\in [t_0,0],\; |\x|\geq c_4 R. 
\end{equation}
We use a partition of unity to construct $W(t,\x)$ so that 
\[
W(t,\x) =\begin{cases} W_1(t,\x),  \quad &|\x|\geq c_4 R+1, \\
-R|\x|+t|\x|^2/2, \quad &|\x|\leq c_4 R. \end{cases}
\]
Clearly $W$ satisfies (2.6) as well.


\subsection{Modified free motion and asymptotic trajectories}

\begin{prop}
Suppose $(x_0,\x_0)$ is backward nontrapping, and let $t_0<0$. 
Then there exists a neighborhood $U$ of $(x_0,\x_0)$ in $\re^{2n}$ such 
that 
\begin{align*}
\x_-(x,\x) &=\lim_{\l\to+\infty} \l^{-1} \y(t_0;x,\l\x), \\
z_-(x,\x) &= \lim_{\l\to+\infty} \bigbra{y(t_0;x,\l\x)
-\pa_\x W(t_0,\y(t_0;x,\l \x))}
\end{align*}
exist for $(x,\x)\in U$. $\x_-(x,\x)$ and $z_-(x,\x)$ are independent of 
$t_0<0$. Moreover, the convergence is uniform in $U$ with its derivatives, and 
\[
S_-:\; (x,\x) \mapsto (z_-,\x_-)
\]
is a local diffeomorphism.
\end{prop}

\begin{rem}
We have already seen $\x_-$ depends only on $(a_{jk}(x))$, and is 
independent of $V(x)$. As we will see in the proof, $z_-$ is also 
independent of $V(x)$, though $W(t,\x)$ does depend on $V(x)$. 
\end{rem}

\begin{proof}
The convergence of $\x_-$ is already shown in Proposition~2.9 and its remark. 
At first, we show 
\[
z^\l(t;x,\x) = y^\l(t;x,\x)-\pa_\x W^\l(t;\y^\l(t;x,\x))
\]
converges as $\l\to\infty$, where 
$W^\l(t,\x)=W(t/\l,\l\x)$ and $t=\l t_0$. 

For $(x,\x)$ near $(x_0,\x_0)$, we choose $\z^\l \in\re^n$ such that 
\[
\y^\l(t;x,\x) =\y^\l(t;-R\z^\l/|\z^\l|,\z^\l), 
\]
and we set
\[
v^\l(s) = y^\l(s;-R\z^\l/|\z^\l|,\z^\l), \quad
w^\l(s) = \y^\l(s;-R\z^\l/|\z^\l|,\z^\l)
\]
for $s\in [t,0]$. Note that $\z^\l$ is a function of $x$, $\x$ and $t=\l t_0$, and 
$\pa_x\pa_\x \z^\l$ is uniformly bounded by virtue of Proposition~2.10 and 
discussion after it. We also set
\begin{align*}
a(s) &= y^\l(s;x,\x) -v^\l(s), \\
b(s) &= \y^\l(s;x,\x)-w^\l(s). 
\end{align*}
We note 
\[
|a(0)|=\bigabs{x+R\z/|\z|} \leq |x|+R, \quad b(t)=0. 
\]
$a$ and $b$ satisfy differential equations:
\begin{align*}
\frac{d}{ds} a(s) &= \frac{\pa p^\l}{\pa\x}(y^\l,\y^\l) 
-\frac{\pa p^\l}{\pa\x}(v^\l,w^\l), \\
\frac{d }{ds} b(s) &= -\biggpare{\frac{\pa p^\l}{\pa x}(y^\l,\y^\l)
-\frac{\pa p^\l}{\pa x}(v^\l,w^\l)},
\end{align*}
where $p^\l(x,\x)= \frac12 \sum_{j,k} a_{jk}(x)\x_j\x_k +\l^{-2} V(x)$. 
Since $\l\geq |s/t_0|$, these imply
\begin{align}
|a'(s)| &\leq c_1 \jap{s}^{-1-\m} |a(s)| +c_1 |b(s)|, \\
|b'(s)| &\leq c_1 \jap{s}^{-2-\m} |a(s)| 
+c_1 \jap{s}^{-1-\m} |b(s)|  
\end{align}
for $s\in[t,0]$ with some $c_1>0$. We note $a(s)=O(\jap{s})$ and $b(s)=O(1)$ 
by Proposition~2.5. Hence by (2.8), we have 
\[
|b(s)| = \biggabs{\int_t^s b'(u)du } \leq c_2 \jap{s}^{-\m}=O(\jap{s}^{-\m}).
\]
Then we substitute this to (2.7) to obtain 
\[
|a(s)| =\biggabs{a(0)-\int^0_s a'(u)du} \leq |x|+R+c_3 \jap{s}^{1-\m}=O(\jap{s}^{1-\m}).
\]
Repeating these, we have $|b(s)|=O(\jap{s}^{-2\m})$ and then $|a(s)|=O(\jap{s}^{1-2\m})$ 
provided $2\m\leq 1$. Iterating this procedure, we arrive at 
$|a(s)|\leq C$ and $|b(s)|\leq C\jap{s}^{-1-\m}$. Moreover, we also have 
\[
|a'(s)|\leq c_4\jap{s}^{-1-\m}.
\]
We recall that $y^\l(s;x,\x) \to \tilde y(s;x,\x)$ as $\l \to \infty$ for each $s$, 
and $\y(t;x,\x)$ converges to $\x_-(x,\x)$ as $\l\to\infty$ since $t=\l t_0$ with $t_0<0$. 
By the uniform continuity of the inverse of $\L(t,\cdot)$, $\z^\l$ converges to $\tilde \z$
as $\l\to\infty$, where 
$\tilde \z$ is given by $\x_-(x,\x)=\x_-(-R\tilde\z/|\tilde\z|,\tilde\z)$.  
Hence, in particular, $v^\l(s)$ converges to $\tilde y(s;-R\tilde\z/|\tilde\z|,\tilde\z)$ 
for each $s$. 
Then by the dominated convergence theorem, we learn
\begin{align*}
&\lim_{\l\to\infty} \bigbra{y(t_0;x,\l\x)-\pa_\x W(t_0,\y(t_0;x,\l\x))} 
= \lim_{\l\to\infty} \bigbra{y^\l (t;x,\x)-v^\l(s)} \\
&= \lim_{\l\to\infty} \biggbrac{x+R\,\frac{\z^\l}{|\z^\l|} -\int_{t}^0\frac{d}{ds} (y^\l(s;x,\x)-v^\l(s))ds}\\
&= x+R\,\frac{\tilde\z}{|\tilde \z|} - \int_{-\infty}^0\frac{d}{ds} (\tilde y(s;x,\x)-\tilde y(s;-R\tilde\z/|\tilde\z|,\tilde\z))ds.
\end{align*}
 Note the right hand side is independent of the potential $V(x)$. 
 
 Next we consider the convergence of the derivatives. As in the proof of Proposition~2.10, for any $\a,\b\in\ze_+^n$, 
 we have 
 \[
| \pa_s\, \pa_x^\a \pa_\x^\b \y^\l(s;x,\x)| \leq C\jap{s}^{-1-\m}, \quad \l t_0\leq s\leq 0. 
\]
Hence, by the dominated convergence theorem, we have 
\begin{align*}
\l^{-1} \pa_x^\a\pa_\x^\b \y(t_0;x,\l\x) &= \pa_x^\a\pa_\x^\b \y^\l(\l t_0;x,\x) \\
&= \x-\int_{\l t_0}^0   \pa_s\, \pa_x^\a \pa_\x^\b \y^\l(s;x,\x) ds \\
&\longrightarrow \ \x-\int^0_{-\infty} \pa_s \pa_x^\a\pa_\x^\b \tilde\y(s;x,\x)ds=\pa_x^\a\pa_\x^\b\x_-(x,\x)
\end{align*}
 as $\l\to\infty$ (cf. Corollary~2.11).
 
 For $z(t;x,\x)$, we prove the convergence by induction. 
 Let $a(s)$ and $b(s)$ as above, and consider $\pa_x^\a\pa_\x^\b a(s)$ and 
 $\pa_x^\a\pa_\x^\b b(s)$. We suppose \[
 \bigabs{\pa_x^\a\pa_\x^\b a(s)}\leq C, \quad  \bigabs{\pa_x^\a\pa_\x^\b b(s)}\leq C\jap{s}^{-1-\m},
 \quad s\in[-\l t_0,0]
\]
for $|\a+\b|<k$ as our induction hypothesis.  
Let $|\a+\b|=k$, and set $A(s)=\pa_x^\a\pa_\x^\b a(s)$ and $B(s)=\pa_x^\a\pa_\x^\b b(s)$. 
Then by inductive computations (from the differential equation for $a(s)$ and $b(s)$), 
we can show (as in the proof of Proposition~2.10), $A(s)$ and $B(s)$ 
satisfy 
\begin{align*}
|A'(s)| &\leq c_1 \jap{s}^{-1-\m} |A(s)| +c_1 |B(s)| +c_1\jap{s}^{-1-\m}, \\
|B'(s)| &\leq c_1 \jap{s}^{-2-\m} |A(s)| +c_1 \jap{s}^{-1-\m} |B(s)| +c_1\jap{s}^{-2-\m} 
\end{align*}
for $s\in[\l t_0,0]$. Note we use a priori estimates: $A(s)=O(\jap{s})$, 
$B(s)=O(1)$, which follow from Proposition~2.10. 
Since $A(0)$ is bounded and $B(\l t_0)=0$, we can use the same  
argument as above (for $a(s)$ and $b(s)$) to conclude  $A(s)=O(1)$ and 
$B(s)=O(\jap{s}^{-1-\m})$, and the induction step is proved. 
Moreover, we have $A'(s)=O(\jap{s}^{-1-\m})$, and the convergence of 
$\pa_x^\a\pa_\x^\b z^\l(t,x,\x) =\pa_x^\a\pa_\x^\b a(t)$ is proved similarly. 

Finally, we prove that $S_-: (x,\x)\mapsto (z_-,\x_-)$ is a local diffeomorphism. 
By the definition, we have 
\[
S_- \exp(TH_p) =S_-
\]
for $T<0$. If $|T|$ is sufficiently large, $\exp(TH_p)$ maps $(x,\x)$ to $(x',\x')$ 
such that $|x'|>>0$ and $x'\cdot\x'<-\d |x'|\,|\x'|$ with some $\d>0$. 
We show $S_-$ is diffeomorphic in a neighborhood of $(x',\x')$ if 
$|x'|$ is sufficiently large. 

We use the above argument for trajectory starting from $(x',\x')$. 
Let $\e>0$ be a small constant which we will specify later. 
Let $0<\m'<\m$. If $|x'|$ is sufficiently large, then $A(s)$ and $B(s)$ above 
(with a new initial condition) satisfy 
\begin{align*}
|A'(s)| &\leq \e c_1 \jap{s}^{-1-\m'} |A(s)| +c_1 |B(s)| +\e c_1\jap{s}^{-1-\m'}, \\
|B'(s)| &\leq \e c_1 \jap{s}^{-2-\m'} |A(s)| +\e c_1 \jap{s}^{-1-\m'} |B(s)| +\e c_1\jap{s}^{-2-\m'} 
\end{align*}
for $s\in[\l t_0,0]$. Then, by carrying out the same argument as above, we learn 
$|A(t)-A(0)|\leq c_2 \e$. In particular, since $z^\l(0)=x+R\z^\l/|\z^\l|$, we have 
\[
|\pa_x (z^\l(t)- x)|\leq c_3 \e, \quad |\pa_\x z^\l(t)| \leq c_3,
\]  
where $t=\l t_0$. We recall, again by Proposition~2.10, we have 
\[
|\pa_x \y^\l(t)|\leq c_3\e, \quad |\pa_\x (\y^\l(t)  - \x)|\leq c_3 \e
\]
if $|x'|$ is sufficiently large. Now if $\e$ is sufficiently small (depending only on $c_3$), 
\[
(x',\x') \mapsto (z^\l(t), \y^\l(t))
\]
has the Jacobian bounded from below by, for example, $1/2$. We now fix $\e>0$, and choose $T$
(and hence $(x', \x')$) accordingly. This Jacobian converges to that of 
$S_-$ as $\l\to\infty$, and hence it is bounded from below by $1/2$. 
Thus we learn that $S_-$ is diffeomorphic in a neighborhood of 
$(x',\x')$ by the inverse function theorem.  
Since $\exp(TH_p)$ is diffeomorphic, this implies $S_-$ is diffeomorphic in a 
neighborhood of $(x,\x)$.  
 \end{proof}
 
 Note the above argument works for the scattering with $\exp tH_k$ also. In fact, the proof is 
 simpler by virtue of the scaling property. For example, $(z^\l(s),\y^\l(s))$ is independent of $\l$,
 the convergence follows immediately from the integrability of the derivative. 


\section{Proof of main theorems}


\subsection{Asymptotic motion and solutions to transport equations}

We denote
\begin{align*}
z(t;x,\x)&= y(t;x,\x)-\pa_\x W(t,\y(t,x,\x)), \\
z^\l(t;x,\x) &= z(t/\l;x,\l \x), \quad \y^\l(t;x,\x)=\y(t/\l;x,\l\x)/\l,
\end{align*}
and also
\begin{align*}
S_t &: \; (x,\x) \mapsto (z(t;x,\x), \y(t;x,\x)), \\
S^\l_t &: \; (x,\x) \mapsto (z^\l(t;x,\x), \y^\l(t;x,\x)).
\end{align*}
$S_t$ (resp. $S^\l_t$) is the Hamilton flow generated by 
\[
\ell(t;x,\x) =p(x+\pa_\x W(t,\x),\x) -\pa_t W(t,\x) 
\]
($\ell^\l(t;x,\x)=\l^{-2} \ell(t/\l,x,\l\x)$, resp.) with the initial condition: 
\[
z(0,x,\x) = x+R\x/|\x|, \quad \y(0;x,\x)=\x
\]
($z^\l(0,x,\x) = x+R\x/|\x|$, $\y^\l(0;x,\x)=\x$, resp.). 
By virtue of the Hamilton-Jacobi equation, we have 
\[
\ell(t,x,\x) =p(x+\pa_\x W(t,\x),\x) -p(\pa_\x W(t,\x),\x)
\]
for sufficiently large $|\x|$. 

Let $f_0(x,\x)$ be a $C_0^\infty$-function supported in a small neighborhood 
of $(x_0+R\x_0/|\x_0|,\x_0)$. We  set 
\[
f_0^\l(x,\x) = f_0(x,\x/\l), 
\]
Then the solution to 
\[
\frac{\pa}{\pa t} f(t;\cdot,\cdot) =-\{\ell,f\}, 
\quad\text{with }\; f(0;x,\x)=f^\l_0(x,\x)
\]
is given by
\[
f(t;x,\x) =f^\l_0\circ S_t^{-1}(x,\x) \qquad \text{for $t\in [t_0,0]$}. 
\]
Similarly, the solution to 
\[
\frac{\pa}{\pa t} f^\l(t;\cdot,\cdot) =-\{\ell^\l,f^\l\},
\quad\text{with }\;  f^\l(0;x,\x)=f_0(x,\x)
\]
is given by 
\[
f^\l(t;x,\x) = f_0\circ (S_t^\l)^{-1}(x,\x)\quad \text{for $t\in [\l t_0,0]$}. 
\]
It is easy to see $f^\l(t;x,\x) = f(\l t;x,\x/\l)$. 
By Proposition~2.12, we learn 
\[
S_-(x,\x) =\lim_{\l\to+\infty} S^\l_{\l t}(x,\x)
\]
exists, and the all the derivatives converges locally uniformly (cf.\ the proof of Proposition~2.12). 
In particular,  we have
\begin{align*}
f_-(x,\x) &= \lim_{\l\to+\infty} f(t;x,\l\x) 
=\lim_{\l\to+\infty} f^\l(\l t;x,\x) \\
&=f_0\circ (S_-)^{-1}(x,\x) \in C_0^\infty(\re^{2n})
\end{align*}
exists and it is independent of $t\in [t_0,0)$. The convergence is locally 
uniform up to its derivatives. 


\subsection{Proof of Theorem~1.2}

At first we consider
\[
v(t) = e^{iW(t,D_x)} e^{-itH}v_0 \quad\text{for } t\in [t_0,0]
\]
with $v_0\in L^2(\re^n)$. $v(t)$ satisfies the evolution equation: 
\begin{align*}
\frac{d}{dt} v(t) &= e^{iW(t,D_x)}\biggbra{ i\frac{\pa W}{\pa t}(t,D_x) -iH}
e^{-itH} v_0 \\
&=-i \biggbra{e^{iW(t,D_x)} H e^{-iW(t;D_x)} -\frac{\pa W}{\pa t}(t,D_x)} v(t).
\end{align*}
Namely, $v(t)$ is a solution to a Schr\"odinger equation with the time-dependent 
Hamiltonian:
\[
L(t) = e^{iW(t,D_x)} H e^{-iW(t,D_x)} -\frac{\pa W}{\pa t}(t,D_x). 
\]
The next lemma is basic in the following analysis. 

\begin{lem}
Let $\n,\rho>0$ and suppose $a\in S\bigpare{\jap{x}^\n\jap{\x}^\rho, 
{dx^2}/{\jap{x}^2}+{d\x^2}/{\jap{\x}^2}}$. Let 
\[
Q = e^{iW(t,D_x)} a(x,D_x) e^{-iW(t,D_x)}. 
\]
Then $Q\in \OPS_K\bigpare{\jap{t\x}^\n\jap{\x}^\rho, 
{dx^2}/{\jap{t\x}^2}+{d\x^2}/{\jap{\x}^2}}$ with any $K\subset\subset \re^n$. 
Let $g(t;x,\x)=\s(Q)$ be the Weyl symbol of $Q$. Then the principal symbol of $Q$
is given by $a(x+\pa_\x W(t,\x),\x)$ and 
\[
g(t;x,\x)-a(x+\pa_\x W(t,\x),\x) \in S_K\biggpare{\jap{t\x}^{\n-2}\jap{\x}^{\rho-2}, 
\frac{dx^2}{\jap{t\x}^2}+\frac{d\x^2}{\jap{\x}^2}}, 
\]
where the remainder is locally bounded in $t$ with respect to the seminorms of 
the symbol class.
\end{lem}

\begin{proof}
The proof is standard pseudodifferential operator calculus, but we sketch it for 
the completeness. Since the Weyl quantization has the same symbol representation 
in the Fourier space as in the configuration space, we may write
\[
\hat A u := \mathcal{F}(a(x,D_x) \check u) =(2\pi)^{-n} \iint e^{-i(\x-\y)\cdot x}
a(x,\tfrac{\x+\y}{2}) u(\y)\,d\y\,dx
\]
for $u\in \mathcal{S}(\re^n)$. By direct computations, we have 
\begin{align*}
& e^{iW(t,\x)} \hat A e^{-iW(t,\x)} u(\x) \\
&\quad= (2\pi)^{-n} \iint e^{i(W(t,\x)-W(t,\y))-i(\x-\y)\cdot x} 
a\bigpare{x,\tfrac{\x+\y}{2}} u(\y) \,d\y\,dx \\
&\quad= (2\pi)^{-n} \iint e^{-i(\x-\y)\cdot (x-\tilde W(t,\x,\y))} 
a\bigpare{x,\tfrac{\x+\y}{2}} u(\y) \,d\y\,dx \\
&\quad = (2\pi)^{-n} \iint e^{-i(\x-\y)\cdot x} 
a\bigpare{x+\tilde W(t,\x,\y),\tfrac{\x+\y}{2}} u(\y) \,d\y\,dx,
\end{align*}
where 
\[
\tilde W(t,\x,\y) =\int_0^1 \pa_\x W(t,s\x+(1-s)\y) \, ds.
\]
We easily see 
\[
\bigabs{\pa_\x^\a \pa_\y^\b \tilde W(t,\x,\y)} \leq C_{\a\b} 
\jap{\tfrac{\x+\y}{2}}^{1-|\a-\b|} \jap{\x-\y}^{1+|\a+\b|}, 
\]
for any $\a,\b\in\ze_+^n$, and $\tilde W(t,\x,\x)=\pa_\x W(t,\x)$. 
Moreover, if $|\a|\geq 2$, by the definition of $W(t,\x)$ and Proposition~2.10, we have 
\begin{equation}
\bigabs{\pa_\x^\a W(t,\x)} \leq C_\a \bigpare{\jap{\x}^{1-|\a|} +|t|\jap{\x}^{2-|\a|}},
\end{equation}
and hence 
\begin{align*}
&\bigabs{\pa_\x^\a \pa_\y^\b \tilde W(t,\x,\y)} \\
&\qquad \leq C_{\a\b}  \biggpare{\bigjap{\tfrac{\x+\y}{2}}^{-|\a+\b|}\jap{\x-\y}^{|\a+\b|} 
+|t|\bigjap{\tfrac{\x+\y}{2}}^{1-|\a+\b|}\jap{\x-\y}^{1+|\a+\b|} }\\
&\qquad \leq C_{\a\b} \bigjap{t\bigpare{\tfrac{\x+\y}{2}}}\bigjap{\tfrac{\x+\y}{2}}^{-|\a+\b|}\jap{\x-\y}^{1+|\a+\b|}.
\end{align*}
We also note 
\[
\bigabs{\tilde W(t,\x,\y)}\geq \d\jap{|t|\tfrac{\x+\y}{2}}
\quad \text{if $\x\cdot\y \geq \frac12|\x|\,|\y|$}
\]
with some $\d>0$. Combining these, we can show 
\begin{align*}
&\bigabs{\pa_x^\a \pa_\x^\b \pa_\y^\c a(x+\tilde W(t,\x,\y),\tfrac{\x+\y}{2})} \\
&\quad \leq C_{\a\b\c} \jap{t(\tfrac{\x+\y}{2})}^{\n-|\a|} 
\jap{\tfrac{\x+\y}{2}}^{\rho-|\b+\c|} \jap{\x-\y}^{|\n|+|\rho|+|\a+\b+\c|}
\end{align*}
for $x\in K\subset\subset \re^n$ and $\x,\y\in \re^n$. 
Then by the asymptotic expansion formula for the simplified symbol, we learn 
that the principal symbol is given by $a(x+\pa_\x W(t,\x),\x)$. Moreover, 
we have 
\[
\bigabs{\pa_x^\a\pa_\x^\b g(t;x,\x)} \leq C_{\a\b} \jap{t\x}^{\n-|\a|}
\jap{\x}^{\rho-|\b|}
\quad \text{for }x\in K, \x\in\re^n,
\]
and the other claims follow from the asymptotic expansion formula. 
\end{proof}

By Lemma~3.1, we learn that the principal symbol of $L(t)$ is given by 
$\ell(t;x,\x)$, and the remainder symbol $r(t;x,\x)$ satisfies
\[
\bigabs{\pa_x^\a \pa_\x^\b r(t;x,\x)} \leq C_{\a\b}
\bigpare{\jap{t\x}^{-\m-2-|\a|}\jap{\x}^{-|\b|} 
+\jap{t\x}^{-\m-|\a|} \jap{\x}^{-2-|\b|}}
\]
for $x\in K\subset\subset \re^n$, $t\in [t_0,0]$. Note that the subprincipal 
symbol vanishes by virtue of the Weyl calculus. 

In order to prove Theorem~1.2, we characterize the wave front set of 
$u_0$ in terms of $u(t_0) =e^{-it_0 H} u_0$ with $t_0<0$. 
Let $a\in C_0^\infty(\re^{2n})$ such that $a(x_0,\x_0)\neq 0$ and supported in 
a small neighborhood of $(x_0,\x_0)$, and set 
\[
a^\l(x,\x) = a(x,\x/\l).
\]
We also set 
\[
A^\l(t) = e^{iW(t,D_x)} e^{-itH} a^\l(x,D_x) e^{itH} e^{-iW(t,D_x)}
\]
for $t\in [t_0,0]$. $A^\l$ satisfies the Heisenberg equation: 
\begin{equation}
\frac{d}{dt} A^\l (t) =-i [L(t), A^\l(t)]. 
\end{equation}
We now construct an asymptotic solution of (3.1) with the initial condition:
\[
A^\l(0) = e^{iW(0,D_x)} a^\l(x,D_x) e^{-iW(0,D_x)} = \tilde a^\l(x,D_x).
\]
We note that the principal symbol of $\tilde a^\l(x,\x)$ is give by 
$a(x-R\hat \x, \x/\l)$, 
and $\tilde a^\l(x,\x)$  is supported in a neighborhood of $(x_0+R\hat \x_0, \l\x_0)$ modulo 
$O(\l^{-\infty})$-terms, where we denote $\hat\x=\x/|\x|$. 

We note that if $A^\l(t)$ is a pseudodifferential operator, the principal 
symbol of the right hand side of (3.2) is given by $-\{\ell,a^\l\}$, where 
$a^\l(t;\cdot,\cdot)$ is the symbol of $A^\l(t)$. Then by the computation in 
Subsection~3.1, we learn that $a^\l\circ S_t^{-1}$ is an approximate solution 
to the transport equation. Actually, we can construct an asymptotic solution 
to (3.2):

\begin{prop}
Let $a\in C_0^\infty(\re^{2n})$ supported in a sufficiently small neighborhood of 
$(x_0,\x_0)$. Then there exists $\g^\l(t;\cdot,\cdot)\in C_0^\infty(\re^{2n})$ such 
that 
\begin{enumerate}
\renewcommand{\theenumi}{\roman{enumi}}
\renewcommand{\labelenumi}{{\rm (\theenumi)}}
\item We write $G^\l(t)=\g^\l(t;x,D_x)$. Then 
\[
G^\l(0)= e^{iW(0,D_x)} a^\l(x,D_x)e^{-iW(0,D_x)} 
\] 
modulo $O(\l^{-\infty})$-terms. 
\item $\g^\l(t;\cdot,\cdot)$ is supported in $S_t[\supp a^\l]$. 
\item For any $\a,\b\in\ze_+^n$, there is $C_{\a\b}>0$ such that 
\[
\bigabs{\pa_x^\a\pa_\x^\b \g^\l(t;x,\x)} \leq C_{\a\b} \l^{-|\b|}, 
\quad t\in[t_0,0],\; x,\x\in\re^n,\; \l\gg 0.
\]
\item The principal symbol of $\g^\l$ is given by $a^\l\circ S_t^{-1}$, i.e., 
\[
\bigabs{\pa_x^\a\pa_\x^\b \bigpare{\g^\l(t;x,\x)-a^\l\circ S_t^{-1}(x,\x)} } 
\leq C_{\a\b} \l^{-1-|\b|}
\]
for $t\in[t_0,0]$,  $x,\x\in\re^n$, $\l\gg 0$.
\item For $t\in [t_0,0]$, 
\[
\biggnorm{\frac{d}{dt} G^\l(t) +i[L(t),G^\l(t)]}_{\mathcal{L}(L^2(\re^n))}
=O(\l^{-\infty}) \quad \text{as $\l\to+\infty$}. 
\] 
\end{enumerate}
\end{prop}

We postpone the proof of Proposition~3.2 to the next subsection, and we 
complete the proof of Theorem~1.2. 

\begin{proof}[Proof of Theorem 1.2]
By Proposition~3.2 and the construction of $L(t)$, we have 
\[
\biggnorm{\frac{d}{dt}\Bigpare{e^{itH} e^{-iW(t,D_x)} G^\l(t) e^{iW(t,D_x)}
e^{-itH} }} \leq C_N \l^{-N}
\]
with any $N$ as $\l\to+\infty$. This implies 
\begin{multline*}
\Bigl\| e^{it_0H} e^{-iW(t_0,D_x)} G^\l(t_0) e^{iW(t_0,D_x)} e^{-it_0H} u_0 \\
-e^{-iW(0,D_x)} G^\l(0) e^{iW(0,D_x)} u_0 \Bigr\| \leq C_N \l^{-N}. 
\end{multline*}
By the condition (i) of Proposition~3.2, we have 
\begin{equation}
\Bigabs{\bignorm{G^\l(t_0) e^{iW(t_0,D_x)} u(t_0)} - \bignorm {a^\l(x,D_x)u_0}}
\leq C_N \l^{-N}, 
\end{equation}
where $u(t)=e^{-itH}u_0$. We note that $\g^\l(t_0;x,\x)$ is supported in 
$S_{t_0}[\supp a^\l]$, and the principal symbol is given by $a^\l\circ S_{t_0}^{-1}$. 
Hence, in particular, 
\begin{equation}
|\g^\l(t_0;x,\x)|\geq \e>0
\end{equation}
for  $|x-z_-(x_0,\x_0)|\leq \d$, $|\x-\l\x_-(x_0,\x_0)|\leq \d\l$ and $\l\gg 0$
with some $\d,\e>0$. 

Now we suppose $(x_0,\x_0)\notin \WF(u_0)$. 
Then by choosing $a$ supported in a sufficiently small neighborhood of $(x_0,\x_0)$, 
we may suppose
\[
\norm{a^\l(x,D_x) u_0} =O(\l^{-\infty}) \quad \text{as $\l\to+\infty$}. 
\]
Then by (3.3) we have 
\begin{equation}
\bignorm{G^\l(t_0) e^{iW(t_0,D_x)} u(t_0)} = O(\l^{-\infty})
\end{equation}
and this implies 
\[
(z_-(x_0,\x_0),\x_-(x_0,\x_0))\notin \WF\bigpare{e^{iW(t_0,D_x)}u(t_0)}
\] 
by virtue of (3.4). 

Conversely, if $(z_-(x_0,\x_0),\x_0(x_0,\x_0))\notin \WF(e^{iW(t_0,D_x)}u(t_0))$ 
then also by taking $a$ supported in a sufficiently small neighborhood of $(x_0,\x_0)$, 
we have (3.5) since $\g^\l(t_0;\cdot,\cdot)$ is supported in $S_{t_0}[\supp a^\l]$ 
modulo $O(\l^{-\infty})$-terms, and it is very close to 
$S_-[\supp a^\l]$ if  $\l$ is large. Then again by (3.3), we have 
$\norm{a^\l(x,D_x) u_0}=O(\l^{-\infty})$, and hence $(x_0,\x_0)\notin \WF(u_0)$. 
\end{proof}


\subsection{Proof of Proposition~3.2}

We note 
\[
R+\d|t\x| \leq |\pa_\x W(t,\x)| \leq R+C|t\x|
\]
for $t\in [t_0,0]$, $\x\in\re^n$ with some $\d,C>0$. 
Using this and (3.1), for any $\a,\b\in\ze_+^n$ and $K\subset\subset\re^n$, 
we have 
\[
\bigabs{\pa_x^\a\pa_\x^\b \ell(t;x,\x)} \leq 
C_{\a\b K}\bigpare{\jap{t\x}^{-1-\m-|\a|}\jap{\x}^{2-|\b|} 
+\jap{t\x}^{1-\m-|\a|}\jap{\x}^{-|\b|}}
\]
for $t\in[t_0,0]$, $x\in K$ and $\x\in\re^n$. 

Let $a_0^\l\in C_0^\infty(\re^{2n})$ such that 
\[
e^{iW(0,D_x)} a^\l(x,D_x) e^{-iW(0,D_x)} = (a_0^\l\circ S_0^{-1})(x,D_x)
\]
modulo $O(\l^{-\infty})$-terms. It is easy to see that the principal symbol 
of $a_0^\l$ is $a^\l(x,\x)$, and that $a_0^\l\in S(1, dx^2+\l^{-2}d\x^2)$. 
We may suppose $\supp a_0^\l =\supp a^\l$. We now set
\[
\g_0(t;x,\x) =a_0^\l\circ S_t^{-1} (x,\x).
\]
Then as we observed in Subsection~3.1, $\g^\l$ satisfies
\[
\frac{\pa }{\pa t} \g_0(t;x,\x)= -\{\ell,\g_0\}(t;x,\x).
\]
We set 
\[
r_0(t;x,\x) = \frac{\pa}{\pa t} \g_0(t;x,D_x)+i[L(t),\g_0(t;x,D_x)].
\]
Then by the asymptotic expansion formula,  
$r_0 \in S(\l^{-1}, dx^2+\l^{-2}d\x^2)$, 
and $r_0$ is supported essentially (i.e., 
modulo $O(\l^{-\infty})$-terms) in $S_t[\supp a^\l]$. 
Next we solve the transport equation: 
\[
\frac{\pa}{\pa t} \g_1(t;x,\x)+\{\ell,\g_1\}(t;x,\x)=-r_0(t;x,\x)
\]
with the initial condition $\g_1(0;x,\x)=0$. It is easy to show that  
$\g_1(t,\cdot,\cdot)\in S(\l^{-1},dx^2+\l^{-2}d\x^2)$ and it is bounded 
in $t\in [t_0,0]$. Moreover, $\g_1$ is supported in $S_t[\supp a^\l]$. 

We set
\[
r_1(t;x,\x) = \frac{\pa}{\pa t} \g_1(t;x,\x) +i[L(t),\g_1(t;x,D_x)] 
+r_0(t;x,D_x), 
\]
then $r_1(t;\cdot,\cdot) \in S(\l^{-2},dx^2+\l^{-2}d\x^2)$ and 
$\supp r_1(t;\cdot,\cdot)\subset S_t[\supp a^\l]$ essentially 
for $t\in[t_0,0]$. We iterate this procedure to obtain 
$\g_j\in S(\l^{-j}, dx^2+\l^{-2}d\x^2)$ such that 
$\supp \g_j(t;\cdot,\cdot)\subset S_t[\supp a^\l]$ essentially 
for $t\in[t_0,0]$.
Then we set
\[
\g^\l(t;x,\x) \sim \sum_{j=0}^\infty \g_j(t;x,\x)
\in S(1,dx^2+\l^{-2}d\x^2)
\]
in the sense of the asymptotic sum as $\l\to+\infty$. By the construction of 
the asymptotic sum, we may suppose 
$\supp \g^\l(t;\cdot,\cdot)\subset S_t[\supp a^\l]$ essentially 
for $t\in[t_0,0]$. 
Now it is straightforward to check $\g$ satisfies the required properties.
\qed


\subsection{Proof of Theorem 1.1}

We denote 
\[
T_t(x,\x) =(x-\pa_\x W(t,\x),\x)
\]
so that 
\[
S_t =T_t \circ \exp tH_p.
\]
We also denote
\[
b_t^\l(x,\x) = a^\l\circ \exp(-tH_p)(x,\x).
\]
Then, in order to prove Theorem~1.1, it suffices to show 
\[
\bignorm{b_{t_0}^\l(x,D_x) u(t_0)} = O(\l^{-\infty}) \quad \text{as $\l\to+\infty$}
\]
if and only if $(x_0,\x_0)\notin \WF(u_0)$, where $u(t)=e^{-itH} u_0$, 
$a$ is supported in a small neighborhood of $(x_0,\x_0)$ and $t_0<0$. 
We note
\[
b_t^\l =a^\l\circ [\exp tH_p]^{-1} = a^\l \circ S_t^{-1}\circ T_t, 
\]
namely, 
\[
b_t^\l(x,\x) = (a^\l\circ S_t^{-1})(x-\pa_\x W(t,\x),\x).
\]

By direct computations as in the proof of Lemma~3.1, we can show 
\[
e^{iW(t,D_x)} b_t^\l(x,D_x) e^{-iW(t,D_x)} = c_t^\l(x,D_x)
\]
where $c_t^\l\in C_0^\infty(\re^n)$ modulo $O(\l^{-\infty})$, and as an 
$h$-pseudodifferential operator (with $h=\l^{-1}$), the principal symbol is given by 
$(a^\l\circ S_t^{-1})(x,\x)=(a\circ (S_{\l t}^\l)^{-1})(x,\x/\l)$. 
Moreover, if we write $\tilde c_t^\l(x,\x) = c_t^\l(x,\l\x)$, 
then $\tilde c_t^\l$ is supported in an arbitrarily small neighborhood of 
$S_-[\supp a]$ if $\l$ is sufficiently large, and the principal symbol is 
$a\circ (S_{\l t}^\l)^{-1}$. We can also show (as in Lemma~3.1) that 
$\tilde c_t^\l$ is bounded in $C_0^\infty(\re^{2n})$ as $\l\to\infty$.
Since 
\[
\bignorm{b_{t_0}^\l(x,D_x) u(t_0)} = \bignorm{\tilde c_t^\l(x,h D_x) e^{iW(t,D_x)} u(t_0)}, 
\]
now Theorem~1.1 follows from Theorem~1.2 combined with the standard 
characterization of the wave front set in terms of $h$-pseudodifferential 
operators. \qed

\appendix
\section{Appendix}

\begin{lem}
Suppose $n\geq 2$, $f\in C^1(\re^n)$ and suppose 
\[
\bigabs{\pa_x f(x)} \leq C\jap{x}^\b, \quad x\in\re^n,
\]
with some $C>0$, $\b\in\re$. Then 
\[
\bigabs{f(x)-f(y)} \leq (\pi/2)C\max(\jap{x}^\b,\jap{y}^\b)|x-y|.
\]
The same estimate holds for $n=1$ if $x\cdot y>0$. 
\end{lem}

\begin{proof}
The claim is obvious if $\b\geq 0$ or $n=1$, and we suppose $n\geq 2$ and 
$\b<0$. Let $|x|\geq |y|\geq 0$ and let 
$S=\bigset{z\in\re^n}{|z|=|y|}$
be the sphere of radius $|y|$ with the center at the origin. 
Let $\ell$ be the (straight) line segment connecting $x$ and $y$.
If $\ell$ and $S$ intersect only at $y$, then we can use the standard 
argument of show
\[
|f(x)-f(y)| =\int_\ell |\nabla f(z)|\,|dz|  \leq C \jap{y}^\b |x-y|.
\]
If $\ell$ and $S$ intersect at $y$ and $y'$, we denote the line segments connecting 
$y$ and $y'$, and $y'$ and $x$, by $\ell'$ and $\ell''$, respectively. 
Then $|x-y|=|\ell|=|\ell'|+|\ell''|$. We note the length of the shortest geodesic 
connecting $y$ and $y'$ on $S$ (which we denote by $\c_1$)  
is equal to or less than $(\pi/2)|\ell'|$.  We set $\c=\c_1+\ell''$, which is a 
piecewise $C^1$-path connecting $y$ and $x$, and 
\[
|\c|\leq (\pi/2)|\ell'|+|\ell''| \leq (\pi/2)|x-y|.
\]
Since $\c$ is contained in $\{z\,|\,|z|\geq |y|\}$, we have 
\[
|f(x)-f(y)| \leq \int_\c |\nabla f(z)|\,  |dz|  \leq C\jap{y}^\b (\pi/2)|x-y|.
\]
\end{proof}


\end{document}